\documentclass[11pt]{article}

\usepackage{amsbsy,amsfonts,amsmath,amssymb,mathrsfs}
\usepackage[T1]{fontenc}
\usepackage{tikz}
\usetikzlibrary{cd,arrows,matrix,positioning}
\usepackage{hyperref}
\usepackage{cleveref}
\usepackage[toc,page,title]{appendix}
\def\itemn#1{\item[\hspace{0.6mm} {\rm (#1)}]}

\renewcommand{\ge}{\geqslant}
\renewcommand{\le}{\leqslant}

\def\toto{\rightrightarrows}
\def\too{\longrightarrow}

\def\into{\hookrightarrow}
\newcommand*{\intoo}{\ensuremath{\lhook\joinrel\relbar\joinrel\rightarrow}}

\def\isomto{\xrightarrow{\,\smash{\raisebox{-0.5ex}{\ensuremath{\scriptstyle\sim}}}\,}}

\renewcommand{\hat}{\widehat}

\makeatletter
\newcommand*{\defeq}{\mathrel{\rlap{%
                     \raisebox{0.3ex}{$\m@th\cdot$}}%
                     \raisebox{-0.3ex}{$\m@th\cdot$}}%
                     =}
\makeatother

\newtheorem{counter}[subsubsection]{$\!\!$}
\newtheorem{subcounter}[subsection]{$\!\!$}

\newenvironment{lemma*}{\begin{subcounter} {\bf Lemma.}}{\end{subcounter}}

\newenvironment{corollary*}{\begin{subcounter} {\bf Corollary.}}{\end{subcounter}}
\newtheorem{theorem}{Theorem}

\newenvironment{remark*}{\begin{subcounter} \rm {\bf Remark.}}{\end{subcounter}}

\newenvironment{notitle*}[1]{\begin{subcounter} {\bf #1.}\rm}{\end{subcounter}}
\newenvironment{proof}{{\flushleft \bf Proof~:}}{\hfill $\square$ \vspace{5mm}}

\DeclareMathOperator{\Spec}{Spec}

\DeclareMathOperator{\colim}{colim}
\DeclareMathOperator{\Hom}{Hom}

\DeclareMathOperator{\Aut}{Aut}
\DeclareMathOperator{\Isom}{Isom}
\DeclareMathOperator{\Mono}{Mono}

\DeclareMathOperator{\id}{id}

\usepackage[toc,page]{appendix}

\DeclareMathOperator{\length}{length}
\DeclareMathOperator{\can}{can}
\DeclareMathOperator{\Sub}{Submt}

 \def\cF{{\cal F}}

\def\cO{{\cal O}}

\def\cZ{{\cal Z}}

\newcommand\CC{\mathbb{C}}

\newcommand\GG{\mathbb{G}}

\newcommand\NN{\mathbb{N}}

\newcommand\ZZ{\mathbb{Z}}

\begin{document}

\begin{center}
{\bf \Large Fixed point stacks under groups
of multiplicative type}

\bigskip

Matthieu Romagny

\bigskip
\today
\end{center}

\bigskip

\begin{center}
\begin{minipage}{17.4cm}
{\small {\bf Abstract.} We prove that if a group scheme of
multiplicative type acts on an algebraic stack with affine,
finitely presented diagonal then the stack of fixed points
is algebraic. For this, we extend two theorems of
\cite{SGA3.2} on functors of subgroups of multiplicative
type, and functors of homomorphisms from a group of
multiplicative type.}
\end{minipage}
\end{center}

\bigskip

We fix a base scheme $S$. For applications in various domains
of algebraic geometry and notably in enumerative geometry,
it is useful to know that an $S$-algebraic stack (in the
sense of Artin) acted on by a group scheme of multiplicative
type has an algebraic stack of fixed points. This was
mentioned in \cite{Ro05}, Remark~3.4 but the case of
general Artin stacks was not considered there. Most recently
Gross, Joyce and Tanaka developed in \cite{GJT20} a framework
for a theory for wall-crossing of enumerative invariants in
$\CC$-linear abelian categories. Provided the abelian category
satisfies certain axioms, these authors explain how a general
mechanism would yield invariants counting semi-stable objects
for weak stability conditions and their wall-crossing formulae.
One of these axioms is for the stack to be algebraic, and a
first choice candidate is the fixed point category of the
compactly supported sheaves on some quasi-projective manifold
(toric surface or Calabi-Yau fourfold) which classifies torus
equivariant compactly supported coherent sheaves.

In the present note we prove such a result of algebraicity
for fixed point stacks under groups of multiplicative type,
answering a question of Arkadij Bojko:

\begin{theorem} \label{th-1}
Let $X\to S$ be an algebraic stack with affine, finitely
presented diagonal. Let $G\to S$ be a finitely presented
group scheme of multiplicative type acting on $X$. Then
the fixed point stack $X^G\to S$ is algebraic, and the
morphism $X^G\to X$ is representable by schemes, separated
and locally of finite presentation.
\end{theorem}

The proof of this soon reduces to the following statement, which
generalizes \cite{SGA3.2}, Exp.~XI,  Th.~4.1 and Cor.~4.2
to non-smooth $H$, thus giving an answer to Exp.~XI, Rem.~4.3.

\begin{theorem} \label{th-2}
Let $H$ be an affine, finitely presented $S$-group scheme.

\smallskip

\noindent {\rm (1)} The functor $\Sub(H)$ of subgroups of
multiplicative
type of $H$ is representable by an $S$-scheme which is
separated and locally of finite presentation. Moreover,
each $S$-quasi-compact closed subscheme of $\Sub(H)$
is affine over $S$.

\smallskip

\noindent {\rm (2)} Let $G$ be a finitely presented $S$-group
scheme of multiplicative type. Then the functor of group
scheme homomorphisms $\Hom(G,H)$ is representable by an
$S$-scheme which is separated and locally of finite
presentation. Moreover, each $S$-quasi-compact closed
subscheme of $\Hom(G,H)$ is affine over $S$.
\end{theorem}

Here are some comments on these results and their proof.
Firstly, it should be emphasized that such statements are very
special to groups of multiplicative type and should not be
expected for more general affine group schemes. They are made
possible by the very strong rigidity properties of groups of
multiplicative type.

Secondly, we observe that Theorem~\ref{th-2} is proved in
\cite{SGA3.2} when $H$ is smooth or more generally when~$H$
is a closed subgroup of a smooth group scheme $H'$ (see
Exp.~XI, Rem.~4.3). This is sufficient to handle
Theorem~\ref{th-1} in the case where $X=[Y/H']$ is a quotient
stack by a smooth group scheme, for then the inertia stack
of $X$ (on which rests the representability of $X^G\to X$,
see Section~\ref{section:algebraicity-of-fixed-pts})
embeds in~$H'$, locally on $X$. However, for the
general case another strategy is needed. What
we do is that we use the density of finite flat
subgroup schemes in groups of multiplicative type to reduce
to a situation where we can apply Grothendieck's theorem on
representability of unramified functors.

We finish this introduction with a remark of a more
historical tone. It has been delightful to read and use
the representability results of \cite{SGA3.2}
and~\cite{Mu65}. This left the author with the impression
to follow a path across the 1960s, seeing in the soil the
seeds of what would ultimately become Artin's criteria for
representability by an algebraic space: the
``auxiliary results
of representability'' of \cite{SGA3.2}, Exp.~XI, \S~3 and
Grothendieck's representability theorem for unramified
functors. May the reader share this pleasure.

The table of contents after the acknowledgements describes
the plan of the article.

\medskip

{\bf Acknowledgements.}
I wish to express warm thanks to Arkadij Bojko, who asked me
whether Theorem~\ref{th-1} held, insisted that I should write
a proof, provided motivation and explanations on the expected
application to computation of enumerative invariants. For
various conversations related to this article, I thank
Alice Bouillet, Pierre-Emmanuel Chaput and Bernard Le Stum.
I acknowledge support from the Centre Henri Lebesgue,
program ANR-11-LABX-0020-01 and would like to thank the
executive and administrative staff of IRMAR and of the
Centre Henri Lebesgue for creating an attractive mathematical
environment.

\tableofcontents

\section{Algebraicity of fixed points stacks}
\label{section:algebraicity-of-fixed-pts}

In this section, we show how to derive Theorem~\ref{th-1}
from Theorem~\ref{th-2}; the rest of the article will be
devoted to the proof of Theorem~\ref{th-2}. Let $X^G$ be
the stack of fixed points as described in \cite{Ro05},
Prop.~2.5. Recall that the sections of $X^G$ over an
$S$-scheme $T$ are the pairs $(x,\{\alpha_g\}_{g\in G(T)})$
composed of an object
$x\in X(T)$ and a collection of isomorphisms
$\alpha_g:x_T\to g^{-1}\,x_T$ satisfying the cocycle
condition $\alpha_{gh}=h^{-1}\alpha_g\circ\alpha_h$:
\[
\begin{tikzcd}[column sep={4em,between origins}]
& h^{-1}\,x_U
\ar[rd,"h^{-1}\alpha_g",start anchor={[xshift=-.5ex]},end anchor={[xshift=-1ex]}] & \\
x_U \ar[ru,"\alpha_h"] \ar[rr,"\alpha_{gh}",swap] & & (gh)^{-1}\,x_U.
\end{tikzcd}
\]
(The $\alpha_g$ here is the $\alpha_{g^{-1}}$ of
\cite{Ro05}; shifting notation simplifies the multiplication
of the group functor~$K$ below.)
It is enough to prove that the morphism $X^G\to X$ is
representable by schemes. For this we fix a point $x:T\to X$
and we study the fibred product
\[
F_G \defeq X^G\times_X T.
\]
This is the functor whose values on a $T$-scheme $U$ are the
collections of isomorphisms $\{\alpha_g\}_{g\in G}$.
In order to prove that $F_G$ is representable by a scheme, we
introduce the functor $K$ defined by
\[
K(U)=\left\{(g,\alpha); \ g\in G(U) \mbox{ and } 
\alpha:x_U\isomto g^{-1}\,x_U \mbox{ an isomorphism}\right\}.
\]
The formula $(g,\alpha) \cdot (h,\beta) \defeq
(gh,h^{-1}\alpha\circ\beta)$ defines a law of multiplication
on $K$:
\[
\begin{tikzcd}[column sep=3em]
x_U \ar[r,"\beta"] & h^{-1}x \ar[r,"h^{-1}\alpha"]
& (gh)^{-1}\,x_U.
\end{tikzcd}
\]
The element $(g,\alpha)=(1,\id_x)$ is neutral for this law.
This makes $K$ a group functor. Moreover, there is a morphism
of group functors $K\to G$, $(g,\alpha)\mapsto g$. This is
representable, affine and finitely presented, because for each
$g:U\to G$ the fibre product $K\times_G U$ is the functor
$\Isom_U(x_U,g^{-1}x_U)$ which is affine and finitely
presented, by assumption on $X$. Since $G$ is also affine
and finitely presented, we see that~$K$ is an affine $T$-group
scheme of finite presentation.

By the definitions, $F_G$ is the functor of group-theoretic
sections of $K\to G$. Because $\Hom(G,G)$ is unramified
and separated (\cite{SGA3.2}, Exp.~VIII, Cor.~1.5), $F_G$
is an open and closed subfunctor of the functor $\Hom(G,K)$
of group scheme homomorphisms. This reduces us to
Theorem~\ref{th-2}, item~(2).

\section{Representability of the functor of homomorphisms:
reductions} \label{section:2}

In Sections~\ref{section:2}
to~\ref{section:representability-using-grothendieck}
we prove item (2) of Theorem~\ref{th-2}: representability
of the functor $\Hom(G,H)$.

\begin{notitle*}{Localization and use of the Density Theorem}
The question of representability is local on $S$ for the
\'etale topology, so we may assume that $S$ is affine
and $G$ is split. Then it is a product $G=N\times \GG_m^r$
where $N$ is finite diagonalizable. But for a product
$G=G_1\times G_2$, the functor $\Hom(G,H)$ is the closed
subfunctor of the product $\Hom(G_1,H)\times \Hom(G_2,H)$
composed of pairs of maps that commute. Using
\cite{SGA3.2}, Exp.~VIII, 6.5.b) we see that it is enough
to consider the factors individually. By \cite{SGA3.2},
Exp.~XI, Prop.~3.12.b) the functor $\Hom(N,H)$ is representable
by an affine $S$-scheme of finite presentation. Hence we
just have to handle the key case $G=\GG_m$.

For a prime number $\ell$ let $S_{\ell}\subset S$ be the open
subscheme where $\ell$ is invertible. We can choose two distinct
primes $\ell,\ell'$ and write $S=S_{\ell}\cup S_{\ell'}$.
Since the question of representability is local on $S$, it is
enough to handle $S_{\ell}$ and $S_{\ell'}$ separately.
In this way we reduce to the case where $\ell\in \cO_S^\times$.

Also, since the functor is locally of finite presentation,
we can reduce to the case where $S$ (still affine) is of finite
type over $\Spec(\ZZ)$.

For each $n>0$, let $\mu_{\ell^n}$ be the kernel of the
endomorphism $\ell^n:\GG_m\to \GG_m$. By \cite{SGA3.2}, Exp.~XI,
Prop.~3.12.b) again, the functor $\Hom(\mu_{\ell^n},H)$ is
representable by an affine $S$-scheme of finite presentation.
By restricting morphisms to the torsion
subschemes, we have a map of functors~:
\[
\varphi:\Hom(\GG_m,H) \too \lim_n \Hom(\mu_{\ell^n},H).
\]
Here the target is representable by an affine scheme, as a
projective limit of affine schemes. Moreover, the Density
Theorem (\cite{SGA3.2}, Exp.~IX, Th\'eor\`eme~4.7 and
Remark~4.10) implies that $\varphi$ is a monomorphism.
We will prove that $\varphi$ representable by using Grothendieck's
theorem on representation of unramified functors. Let $T$ be
an $S$-scheme and let $\{u_n:\mu_{\ell^n,T}\to H_T\}$ be a collection
of morphisms of $T$-group schemes. We want to prove that the
fibred product
\[
\Hom(\GG_m,H) \times_{\lim_n \Hom(\mu_{\ell^n},H)} T
\]
is representable. For this we change our notation, renaming
$T$ as $S$, and we end up with the following setting.
\end{notitle*}

\begin{notitle*}{Summary of situation} \label{summary}
After these reductions, we start anew with:
\begin{itemize}
\item a prime number $\ell$,
\item a $\ZZ[1/\ell]$-algebra of finite type $R$
and $S=\Spec(R)$,
\item a finitely presented affine $S$-group scheme $H$,
\item a family of compatible morphisms of $S$-group
schemes $\{u_n:\mu_{\ell^n}\to H\}$.
\end{itemize}
The compatibility of the $u_n$'s means that $u_{n+1}$
extends $u_n$ for each $n$. We study the functor
\[
F(T)=\{\mbox{morphisms } f:\GG_m\to H
\mbox{ that extend the $u_n$, $n\ge 0$}\}.
\]
The Density Theorem implies that $F\to S$ is a monomorphism of
functors. That is, in fact $F(T)=\{f\}$ is the one-element set
composed of the unique morphism $f:\GG_{m,T}\to H_T$ that extends
the $u_n$ if there is one, and $F(T)=\varnothing$ otherwise.
In particular $F\to S$ is formally unramified. In order to prove
the representability of $F$ in
Section~\ref{section:representability-using-grothendieck},
we will use Grothendieck's theorem on unramified functors,
verifying the eight relevant axioms. For this we need some
preparations which we do in the next section.
\end{notitle*}

\section{Descent along schematically dominant morphisms}
\label{section:3}

Again let $F$ be the functor of group homomorphisms
$f:\GG_m\to H$ extending the maps $u_n$. The verification
of Grothendieck's axioms for $F$ will be based
to a large extent on the following fact: the map
$F(T)\to F(T')$ is an isomorphism for all schematically
dominant morphisms of schemes $T'\to T$. This is
Lemma~\ref{lemma:descent-by-sch-dominant} below. Its proof
will use a variation on the argument used to show that
formal homomorphisms from a group scheme of multiplicative type
to an affine group scheme are algebraic, see \cite{SGA3.2}
Exp.~IX, \S~7. It is the purpose of this subsection to
settle this.

We work over a $\ZZ[1/\ell]$-algebra $A$.
The argument in {\em loc. cit.} uses the natural
embedding of the Hopf algebra $A[\GG_m]$, which as a module is
the direct sum $A^{(\ZZ)}$, into the direct product $A^{\ZZ}$.
In our context, the useful information we have comes from the
embedding of $A[\GG_m]$ into the limit of the Hopf algebras of
the $\mu_{\ell^n}$'s. We need to relate the two embeddings.

\begin{notitle*}{The ind-scheme of $\ell$-power roots of unity}
We need some facts on the ind-scheme
of $\ell$-power roots of unity:
\[
\mu_{\ell^\infty}=\colim \mu_{\ell^n}.
\]
We consider its function algebra
$A[\mu_{\ell^\infty}]=\lim A[\mu_{\ell^n}]$ and
the canonical injective morphism:
\[
c:A[\GG_m]\intoo A[\mu_{\ell^\infty}].
\]
Contemplated from the right
angle, this map is naturally isomorphic to the inclusion of
$A$-modules $A^{(\ZZ)}\into A^{\ZZ}$ of the direct sum into
the direct product. For later use, we need to make this precise.
For each $n\ge 0$ let $\mu_{\ell^n}^*\subset \mu_{\ell^n}$ be the
subscheme of primitive roots of unity; we have
$\mu_{\ell^n}^*=\Spec(A[z]/(\Phi_n))$ where $\Phi_n$
is the cyclotomic polynomial.  Since $\ell$
is invertible in the base ring, the $\Phi_n$ are pairwise
strongly coprime in $A[z]$ and the factorization
$z^{\ell^n}-1 =\Phi_0\cdots\Phi_n$ gives rise to isomorphisms
of algebras:
\[
A[\mu_{\ell^n}]=\prod_{0\le i\le n} \frac{A[z]}{(\Phi_i)}
\ ,\qquad
A[\mu_{\ell^\infty}]=\prod_{i\in\NN} \frac{A[z]}{(\Phi_i)}.
\]
\end{notitle*}

\begin{notitle*}{$\Phi$-adic expansions}
Let $n\ge 0$. Let
$I_n:=\{-\lfloor \varphi(\ell^n)/2\rfloor,\dots,\varphi(\ell^n)-\lfloor \varphi(\ell^n)/2\rfloor-1\}$ where $\varphi$ is Euler's
totient function. The $A$-module $C_n=\oplus_{i\in I_n} Az^i$
is finite free of rank $\varphi(\ell^n)=\deg(\Phi_n)$. It follows
that the composition $C_n\into A[z^{\pm 1}]\to A[z]/(\Phi_n)$ is
an isomorphism of $A$-modules, and we can choose~$C_n$ as a module
of representatives of residue classes for Euclidean division
mod $\Phi_n$. As a result, for each Laurent polynomial
$P\in A[z^{\pm 1}]$ there are unique $q\in A[z^{\pm 1}]$ and
$r\in C_n$ such that $P=\Phi_n q+r$. By writing such divisions~:
\begin{align*}
P=\Phi_0q_0+r_0 \quad(r_0\in C_0)\\
q_0=\Phi_1q_1+r_1 \quad(r_1\in C_1)\\
q_1=\Phi_2q_2+r_2 \quad(r_2\in C_2)
\end{align*}
and so on, we reach a $\Phi$-adic expansion
$P=r_0+r_1(z-1)+\dots +r_{d+1}(z^{\ell^d}-1)$ for some $d$.
This gives a direct sum decomposition by the submodules $D_0=A$ and
$D_n=C_n\cdot (z^{\ell^{n-1}}-1)$ for all $n\ge 1$:
\[
A[z^{\pm 1}]=\bigoplus_{n\ge 0} D_n.
\]
\end{notitle*}

\begin{notitle*}{The map $\theta_A:A^{\ZZ}\to A[\mu_{\ell^\infty}]$}
As a module, the ring $A[z^{\pm 1}]$ is the direct sum
$A^{(\ZZ)}$; we embed it into the direct product $A^{\ZZ}$.
Compatibly with the expression $A^{(\ZZ)}=\oplus_{n\ge 0} D_n$,
we have a direct product decomposition
$A^{\ZZ}=\prod_{n\ge 0} D_n$. Also we view the algebra
$A[\mu_{\ell^\infty}]$ as the product
$\prod_{i\in\NN} \frac{A[z]}{(\Phi_i)}$. The restriction
$D_n \into A[z^{\pm 1}] \to A[z]/(\Phi_n)$ of the residue
class map is an isomorphism; we abuse notation slightly by
omitting it from the notation in defining the following map:
\[
\begin{tikzcd}[column sep=25,
    ,row sep = 0ex
    ,/tikz/column 1/.append style={anchor=base east}
    ,/tikz/column 2/.append style={anchor=base west}
    ]
A^{\ZZ}=\underset{n\in \NN}{\prod} D_n \ar[r,"\theta_A"] &
A[\mu_{\ell^\infty}]=\underset{n\in \NN}{\prod}\frac{A[z]}{(\Phi_n)} \\
(P_0,P_1,P_2,\dots) \ar[r,mapsto] & (P_0,P_0+P_1,P_0+P_1+P_2,\dots).
\end{tikzcd}
\]
\end{notitle*}

\begin{lemma*} \label{lemma:map-theta}
{\rm (1)} The map $\theta_A$ is an isomorphism of $A$-modules
and fits in a commutative diagram:
\[
\begin{tikzcd}[column sep=20,row sep=30]
A^{\ZZ} \ar[rr,"\theta_A","\sim"'] & &
A[\mu_{\ell^\infty}] \\
& A^{(\ZZ)}
\ar[lu,hook',"{\can}",start anchor={[xshift=-3ex]},end anchor={[xshift=-1ex]}]
\!=\!A[\GG_m] \ar[ru,hook,"c"',start anchor={[xshift=3ex]},end anchor={[xshift=2ex]}]
&
\end{tikzcd}
\]
Here $\can$ is the canonical inclusion of the direct
sum into the product.

\smallskip

\noindent {\rm (2)} The isomorphism $\theta_A$ is
functorial in $A$, that is, for each morphism of
$\ZZ[1/\ell]$-algebras $A\to B$ there is a cartesian
diagram compatible with the maps $\can_A,\can_B$ and $c_A,c_B$:
\[
\begin{tikzcd}[column sep=30]
A[\mu_{\ell^\infty}] \ar[d] \ar[r,"\theta_A"]
& A^{\ZZ} \ar[d] \\
B[\mu_{\ell^\infty}] \ar[r,"\theta_B"] & B^{\ZZ}
\end{tikzcd}
\]
In particular, if $A\to B$ is injective and $P\in B^{(\ZZ)}$,
we have $\can(P) \in A^{\ZZ}$ if and only if
$c(P)\in A[\mu_{\ell^\infty}]$.
\end{lemma*}

\begin{proof}
(1) Let $P\in A[z^{\pm 1}]$ be written in the form
$P=P_0+P_1+\dots+P_d$ with $P_i\in D_i$ for all~$i$. The canonical
injection $\can:\oplus D_n \into \prod D_n$ sends $P$ to
$(P_0,P_1,P_2,\dots)$ 
while the map $c$ sends $P$ to $(P_0,P_0+P_1,P_0+P_1+P_2,\dots)$.
In other words, $c(P)=\theta(\can(P))$ as desired.
The fact that $\theta$ is an $A$-linear isomorphism follows
from the fact that it is defined by a triangular unipotent matrix.

\smallskip

\noindent {\rm (2)} The commutative diagram exists by the
very construction; it is cartesian
because $\theta_A$ and $\theta_B$ are isomorphisms.
\end{proof}

We can now prove that the objects of the functor $F$ descend
along schematically dominant morphisms. For the latter notion,
we refer the reader to \cite{EGA} IV$_3$.11.10.

\begin{lemma*} \label{lemma:descent-by-sch-dominant}
Let $T'\to T$ be a morphism of $S$-schemes and assume that either:
\begin{trivlist}
\itemn{1} $T'\to T$ is schematically dominant, or
\itemn{2} $T=\Spec(A)$ is affine, $T'=\amalg_i \Spec(A_i)$
is a disjoint sum of affines, with $A\to \prod_i A_i$ injective.
\end{trivlist}
Then the map $F(T)\to F(T')$ is bijective.
\end{lemma*}

The result is easier when $T'\to T$ is quasi-compact,
but the general case will be crucial for us.

\begin{proof}
(1) Since $F(T)$ has at most one point, the map $F(T)\to F(T')$ is
injective and it is enough to prove that it is surjective.
By fpqc descent of morphisms, this holds when $T'\to T$
is a covering for the fpqc topology. Applying this remark
with chosen Zariski covers $\amalg T_i\to T$ and
$\amalg_{i,j} T'_{ij}\to \amalg_i T_i\times_T T'\to T'$,
we see that the vertical
maps in the following commutative square are bijective:
\[
\begin{tikzcd}
F(T) \ar[r]
\ar[d,"\sim"{above,label distance=1pt,sloped,pos=.4}] & F(T') \ar[d,"\sim"{above,label distance=1pt,sloped,pos=.4}] \\
\prod F(T_i) \ar[r] & \prod F(T'_{ij}).
\end{tikzcd}
\]
Choosing $T=\Spec(A)$ and $T'_{ij}=\Spec(A'_{ij})$ afffine,
the assumption that $T'\to T$ is schematically dominant implies
that $A\to \prod A'_{ij}$ is injective. This way we reduce to case (2).

\smallskip

\noindent (2)
We start with an element of $F(T')$, i.e. a family of morphisms
of $A_i$-group schemes $f_i:\GG_{m,A_i}\to H_{A_i}$ each of which
extends the morphisms $u_n:\mu_{n,A_i}\to H_{A_i}$, $n\ge 0$.
For simplicity, in the sequel we write again
$f_i:\cO_H\otimes A_i\to A_i[z^{\pm 1}]$ and
$u_n:\cO_H\to R[z]/(z^{\ell^n}-1)$ the corresponding comorphisms
of Hopf algebras; this should not cause confusion. For each
$R$-algebra $A$ we also write
$u_{\infty,A}:\cO_H\otimes A\to A[\mu_{\ell^{\infty}}]$ for
the product of the $u_{n,A}$. These fit in a commutative diagram:
\[
\begin{tikzcd}[column sep=30]
\cO_H\otimes A_i \ar[d,"f_i"{swap}] \ar[rd,"u_{\infty,A_i}"]
& & \\
A_i[z^{\pm 1}] \ar[r,"c_{A_i}"{swap}] &
A_i[\mu_{\ell^\infty}].
\end{tikzcd}
\]
We now reduce to the case where $A$ and the $A_i$ are
noetherian. For this let $L$ resp. $L_i$ be the image of
$R\to A$, resp. of $R\to A_i$. Being quotients of $R$,
the rings $L$ and $L_i$ are noetherian. Moreover, since
$A\to \prod A_i$ is injective then so is $L\to \prod L_i$.
Since the $u_n$ are defined over $R$ hence over $L$,
we have a commutative diagram:
\[
\begin{tikzcd}[column sep=30,row sep=30]
\cO_H\otimes L_i \ar[d] \ar[rrd,"u_{\infty,L_i}"]
\ar[rd,"{f_{i,L_i}}"{swap} label distance=1mm,dotted] & \\
\cO_H\otimes L_i \ar[rd,"f_i"{swap}] &
L_i[z^{\pm 1}] \ar[r,"c_{L_i}"{swap}] \ar[d] &
L_i[\mu_{\ell^\infty}] \ar[d] \\
& A_i[z^{\pm 1}] \ar[r,"c_{A_i}"{swap}] &
A_i[\mu_{\ell^\infty}].
\end{tikzcd}
\]
Since the lower right square is cartesian, there is an induced
dotted arrow. In this way we see that $f_i$ is actually defined
over $L_i$. So replacing $A$ (resp. $A_i$) by $L$ (resp.
$L_i$), we obtain the desired reduction.

Let $\hat A\defeq \prod_i A_i$. Taking products over $i$,
we build a commutative diagram:
\[
\begin{tikzcd}[column sep=40,row sep=10]
\cO_H\otimes A \ar[d] \ar[rr,"u_{\infty,A}"]
\ar[rdd,dotted] &
& A[\mu_{\ell^\infty}] \ar[dd,hook] \\
\cO_H\otimes \hat A \ar[d] & & \\
\prod_i(\cO_H\otimes A_i) \ar[r,"\prod f_i"]
& \prod_i(A_i[z^{\pm 1}]) \ar[r,hook,"{\prod c_{A_i}}"]
& \hat A[\mu_{\ell^\infty}]
\end{tikzcd}
\]
Henceforth we set $C\defeq \cO_H\otimes A$ and we write
$\Phi_0$ the dotted composition in the diagram above.
What the diagram shows is that
\[
\begin{tikzcd}[column sep=45]
C \ar[r,"\Phi_0"] &
\prod_i(A_i[z^{\pm 1}]) \ar[r,hook,"{\prod c_{A_i}}"]
& \hat A[\mu_{\ell^\infty}]
\end{tikzcd}
\]
factors through $A[\mu_{\ell^\infty}]$. According to
Lemma~\ref{lemma:map-theta}(2) applied with
$B=\hat A$, this implies that
\[
\begin{tikzcd}[column sep=45]
C \ar[r,"\Phi_0"] &
\prod_i(A_i[z^{\pm 1}]) \ar[r,hook,"{\prod \can_{A_i}}"]
& {\hat A}{}^{\,\ZZ}
\end{tikzcd}
\]
factors through $A^{\ZZ}$, providing a map
$\Phi:C \to A^{\ZZ}$.
From the diagrams expressing the fact that the~$f_i$
respect the comultiplications, taking products over $i$,
we obtain a commutative diagram:
\[
\begin{tikzcd}
C \ar[rr,"\Phi"] \ar[d] & & A^{\ZZ} \ar[d] \\
C\otimes_A C \ar[r,"\Phi\otimes\Phi"] &
A^{\ZZ}\otimes_A A^{\ZZ} \ar[r] & A^{\ZZ\times\ZZ}.
\end{tikzcd}
\]
Let $g\in C$ and write $\Phi(g)=(a_m)_{m\in\ZZ}$. Since
$A$ is noetherian, Lemme~7.2 of \cite{SGA3.2}, Exp.~IX
is applicable and shows that only finitely many of the
$a_m$ are nonzero, that is $\Phi(g)\in A[z^{\pm 1}]$.
Therefore $\Phi$ gives rise to a map $f:C\to A[z^{\pm 1}]$.
The fact that~$f$ respects the comultiplication of the
Hopf algebras follows immediately by embedding
$A[z^{\pm 1}]\otimes A[z^{\pm 1}]$ into
${\hat A}[z^{\pm 1}]\otimes {\hat A}[z^{\pm 1}]$ where the
required commutativity holds by assumption.
The fact that $f$ respects the counits is equally clear.
\end{proof}

\section{Representability of the functor of homomorphisms:
proof}
\label{section:representability-using-grothendieck}

We come back to the setting of~\ref{summary} and we proceed
to show that $F\to S$ is representable; this will complete
the proof of item (2) of Theorem~\ref{th-2}. For this, we
use Grothendieck's theorem on representation of unramified
functors, which we begin by recalling.

\begin{notitle*}{Representation of unramified functors:
statement}
In what follows, all affine schemes $\Spec(A)$ that appear
are assumed to be $S$-schemes, and we write $F(A)$ instead
of $F(\Spec(A))$.

\begin{theorem}{\bf (Grothendieck~\cite{Mu65})}
Let $S$ be a locally noetherian scheme and $F$ a set-valued
contravariant functor on the category of $S$-schemes. Then $F$
is representable by an $S$-scheme which is locally of finite
type, unramified and separated if and only if Conditions
$(F_1)$ to $(F_8)$ below hold.
\begin{itemize}
\item[$(F_1)$] The functor $F$ is a sheaf for the fpqc topology.
\item[$(F_2)$] The functor $F$ is locally of finite presentation;
that is,
for all filtering colimits of rings $A=\colim A_\alpha$,
the map $\colim F(A_\alpha)\to F(A)$ is bijective.
\item[$(F_3)$] The functor $F$ is effective; that is, for all noetherian
complete local rings $(A,m)$, the map $F(A)\to \lim F(A/m^k)$ is
bijective.
\item[$(F_4)$] The functor $F$ is homogeneous; more precisely, for all exact
sequences of rings $A\to A'\toto A'\otimes_AA'$
with $A$ local artinian, $\length_A(A'/A)=1$ and trivial residue
field extension $k_A=k_{A'}$, the diagram
$F(A)\to F(A')\toto F(A'\otimes_AA')$ is exact.
\item[$(F_5)$] The functor $F$ is formally unramified.
\item[$(F_6)$] The functor $F$ is separated; that is, it satisfies the valuative
criterion of separation.
\end{itemize}
For the last two conditions we let $A$ be a noetherian ring,
$N$ its nilradical, $I$ a nilpotent ideal such that $IN=0$,
$T=\Spec(A)$, $T'=\Spec(A/I)$. We assume that $T$
is irreducible and we call $t$ its generic point.
\begin{itemize}
\item[$(F_7)$] Assume moreover that $A$ is complete one-dimensional local
with a unique associated prime. Then any point $\xi':\Spec(A/I)\to F$ such
that
\[
\xi'_t:\Spec((A/I)_t)\too \Spec(A/I)\too F
\]
can be lifted to a point $\xi^*:\Spec(A_t)\to F$, can be lifted
to a point $\xi:\Spec(A)\to F$.
\item[$(F_8)$] Assume that $\xi':\Spec(A/I)\to F$ is such that
\[
\xi'_t:\Spec((A/I)_t)\too \Spec(A/I)\too F
\]
can not be lifted to any subscheme of $\Spec(A_t)$ which is strictly
larger than $\Spec((A/I)_t)$. Then there exists a nonempty open set
$W\subset T$ such that for all open subschemes $W_1\subset T$ contained
in~$W$, the restriction $\xi'_{|W'_1}:W'_1\to F$ (with $W'_1=W_1\times_TT'$)
can not be lifted to any subscheme of~$W_1$ which is strictly larger than
$W'_1$.
\end{itemize}
\end{theorem}

We now start to verify the conditions one by one.
\end{notitle*}

\begin{notitle*}{Conditions $(F_1)$, $(F_4)$, $(F_5)$, $(F_6)$, $(F_7)$} We begin by checking the conditions which turn out
easy in our case.

\medskip

\noindent $(F_1)$ This follows from fpqc descent, see e.g.
\cite{SGA1}, Exp.~VIII, Th.~5.2.

\smallskip

\noindent $(F_4)$ Since $A\to A'$ is injective,
by Lemma~\ref{lemma:descent-by-sch-dominant}
the map $F(A)\to F(A')$ is a bijection. This gives a statement
which is much stronger than the $(F_4)$ in the theorem.

\smallskip

\noindent $(F_5)$ Since $F\to S$ is a monomorphism, it is
formally unramified.

\smallskip

\noindent $(F_6)$ Since $F\to S$ is a monomorphism, it is separated.

\smallskip

\noindent $(F_7)$ Since $A$ has a unique associated prime, the map
$\Spec(A_t)\to\Spec(A)$ is schematically dominant.
Hence by Lemma~\ref{lemma:descent-by-sch-dominant}, the point
$\xi^*:\Spec(A_t)\to F$ automatically extends to $\Spec(A)$.
\end{notitle*}

\begin{notitle*}{Condition $(F_2)$}
Let $\cF:=\Hom(\GG_m,H)$ be the functor of {\em all} morphisms
of group schemes $\GG_m\to H$ --- that is, not just those that
extend the collection $u_n$. It is standard that $\cF$ is
locally of finite presentation (see \cite{EGA} IV$_3$.8.8.3).
Should the affine scheme $\lim_n \Hom(\mu_{\ell^n},H)$ be
locally of finite type over $S$, it would follow that $F\to S$
is locally of finite presentation (\cite{EGA} IV$_1$.1.4.3(v)).
However this is not the case in general, and the verification
of $(F_2)$ needs more work.

So let $A=\colim A_\alpha$ be a filtering colimits of rings.
We want to prove that $\colim F(A_\alpha)\to F(A)$ is bijective.
We look at the diagram
\[
\begin{tikzcd}
\colim F(A_\alpha) \ar[r] \ar[d,hook] & F(A) \ar[d,hook] \\
\colim \cF(A_\alpha) \ar[r,"\sim"] & \cF(A).
\end{tikzcd}
\]
Since $\cF$ is locally of finite presentation, the bottom
row is an isomorphism. We deduce that the upper row is
injective. We shall now prove that the upper row is
surjective, and in fact that the diagram is cartesian.

\begin{lemma*}
Let $f:\GG_{m,A}\to H_A$ be a morphism extending
$u_{n,A}:\mu_{\ell^n,A}\to H_A$
for all $n\ge 0$. Then there exists an index $\alpha$ such
that $f$ descends to a map $f_\alpha:\GG_{m,A_\alpha}\to H_{A_\alpha}$
extending $u_{n,A_\alpha}$ for all $n\ge 0$.
\end{lemma*}

\begin{proof}
Since $G$ and $H$ are finitely presented, the morphism $f$
is defined at finite level, that is there exists an index $\alpha$
and a morphism of $A_\alpha$-group schemes
$g:\GG_{m,A_\alpha}\to H_{A_\alpha}$
whose pullback along $\Spec(A)\to\Spec(A_\alpha)$ is $f$.
Since the groups are affine, the morphism $g$ is given by a map of
rings $g^\sharp:\cO_H\otimes A_\alpha\to A_\alpha[z^{\pm 1}]$. Fix
a presentation $\cO_H=R[x_1,\dots,x_s]/(P_1,\dots,P_t)$. Then~:
\begin{itemize}
\item $u_n^\sharp$ is determined by the
elements $z_{n,j}\defeq u_n^\sharp(x_j)\in R[z]/(z^{\ell^n}-1)$
satisfying $P_k(z_{n,1},\dots,z_{n,s})=0$ for $k=1,\dots,t$,
\item $g^\sharp$ is determined by the elements
$y_j=g^\sharp(x_j)\in A_\alpha[z^{\pm 1}]$ satisfying
$P_k(y_1,\dots,y_s)=0$ for $k=1,\dots,t$.
\end{itemize}
Moreover, saying that $f$ extends $u_{n,A}$ is just saying that
$z_{n,j}=\pi_n(y_j)$ in $A[z]/(z^{\ell^n}-1)$, for all $j$,
where
\[
\pi_n:A_\alpha[z^{\pm 1}]\to A[z]/(z^{\ell^n}-1)
\]
is the projection. So we have to prove that we may enlarge
the index $\alpha$ in such a way that $g$ extends
$u_{n,A_\alpha}$ for all $n\ge 0$.

For $n_0\ge 0$ an integer, let
$J_0:=\{-\lfloor\ell^{n_0}/2\rfloor,\dots,\ell^{n_0}-\lfloor\ell^{n_0}/2\rfloor-1\}$
and $E_0=\oplus_{i\in J_0} R z^i$. Then $\pi_{n_0|E_0}$
is an isomorphism and for all $n\ge n_0$ we can define
\[
\chi_n=\pi_n\circ (\pi_{n_0|E_0})^{-1}:
R[z]/(z^{\ell^{n_0}}-1) \too R[z]/(z^{\ell^n}-1).
\]
By base change, these objects are defined over any $R$-algebra.
We choose $n_0$ large enough so that $E_0\otimes_R A_\alpha$
contains the Laurent polynomials $y_1,\dots,y_s$.

In the present context, the condition that $f$ extends all the
maps $u_{n,A}:G_{n,A}\to H_A$ is a finiteness constraint imposed
by $f$ on $\{u_n\}$ (whereas in other places of our arguments
it is best seen as a condition imposed by $\{u_n\}$ on $f$).
Indeed, from the relations $z_{n,j}=\pi_n(y_j)$ in $A$, we deduce
that
\[
z_{n,j}=\chi_n(z_{n_0,j}) \quad\mbox{in }
A[z]/(z^{\ell^n}-1) \quad\mbox{for all } n\ge n_0,
\]
namely
$\chi_n(z_{n_0,j})=(\pi_n\circ (\pi_{n_0|E_0})^{-1})(\pi_{n_0}(y_j))
=\pi_n(y_j)=z_{n,j}$. We claim that we may increase $\alpha$
to achieve that these equalities hold in $A_\alpha[z]/(z^{\ell^n}-1)$,
for all $n\ge n_0$ and all $j$. In order to see this, note that
the elements $\delta_{n,j}:=z_{n,j}-\chi_n(z_{n_0,j})$ are defined
over $R$, and as we have just proved, they belong to the kernel of
the morphism $R[z]/(z^{\ell^n}-1)\to A[z]/(z^{\ell^n}-1)$. Let
$I\subset R$ be the ideal generated by the coefficients of the
expressions of $\delta_{n,j}$ on the monomial basis, for varying
$n\ge n_0$ and $j$. Since $R$ is noetherian, $I$ is generated by
finitely many elements. These elements vanish in $A$, hence they
vanish in~$A_\alpha$ provided we increase $\alpha$ a little, whence
our claim.

\smallskip

The relations $z_{n,j}=\pi_n(y_j)$ in $A[z]/(z^{\ell^n}-1)$
with $j=1,\dots,s$ and $n\le n_0$ being finite in number, we may
increase $\alpha$ so as to ensure that all of them hold in
$A_{\alpha}[z]/(z^{\ell^{n}}-1)$.
Then for $n\ge n_0$ we have
\[
z_{n,j}=\chi_n(z_{n_0,j})=\chi_n(\pi_{n_0}(y_j))=\pi_n(y_j)
\quad\mbox{in } A_{\alpha}[z]/(z^{\ell^{n_0}}-1)
\]
again. That is, $g$ extends
the maps $u_{n,A_{\alpha}}$ for all $n\ge 0$.
\end{proof}
\end{notitle*}

\begin{notitle*}{Condition $(F_3)$}
Let $(A,m)$ be a noetherian complete local ring. We want to prove
that the map $F(A)\to \lim F(A/m^k)$ is bijective. We write again
$\cF:=\Hom(\GG_m,H)$. We look at the diagram
\[
\begin{tikzcd}
F(A) \ar[r] \ar[d,hook] & \lim F(A/m^k) \ar[d,hook] \\
\cF(A) \ar[r,"\sim"] & \lim \cF(A/m^k).
\end{tikzcd}
\]
From \cite{SGA3.2}, Exp.~IX, Th.~7.1 we know that $\cF$ is effective,
that is the bottom arrow is bijective. We deduce that
the upper row is injective. We shall now prove that the upper row is
surjective, and in fact that the diagram is cartesian. So let
$f_k:\GG_{m,A/m^k}\to H_{A/m^k}$ be a collection of $A/m^k$-morphisms
such that $f_k$ extends $u_{n,A/m^k}:\mu_{\ell^n,A/m^k}\to H_{A/m^k}$
for all $n\ge 0$, and let $f:\GG_{m,A}\to H_A$ be a morphism that
algebraizes the $f_k$. We must prove that $f$ extends $u_{n,A}$,
for each $n$. For this let $i_n:\mu_{\ell^n,A}\to\GG_{m,A}$ be
the closed immersion. The two maps $f\circ i_n$ and $u_n$
coincide modulo $m^k$ for each $k\ge 1$, hence so do the
morphisms of Hopf algebras
\[
(f\circ i_n)^\sharp,u_n^\sharp:
\cO_H\otimes A\to A[z]/(z^{\ell^n}-1).
\]
Since $A[z]/(z^{\ell^n}-1)$ is separated for the $m$-adic
topology, we deduce that $(f\circ i_n)^\sharp=u_n^\sharp$
and hence $f\circ i_n=u_n$. This concludes the argument.

\begin{remark*}
We could also appeal to the following more general result
extending the injectivity part of \cite{EGA}~III$_1$.5.4.1:
{\em let $(A,m)$ be a
noetherian complete local ring, and $S=\Spec(A)$. Let $X,Y$
be $S$-schemes of finite type with $X$ pure and $Y$ separated.
Let $f,g:X\to Y$ be $S$-morphisms. If we have the equality of
completions $\hat f=\hat g$, then $f=g$.} For the notion of
a pure morphism of schemes we refer to Raynaud and
Gruson~\cite{RG71} and Romagny~\cite{Ro12}. As to the proof
of the italicized statement,
by \cite{EGA1new}, 10.9.4 the morphisms~$f$ and $g$ agree in
an open neighbourhood of $\Spec(A/m)$. Then the arguments in
the proof of \cite{Ro12}, Lemma~2.1.9 apply verbatim.
\end{remark*}
\end{notitle*}

\begin{notitle*}{Condition $(F_8)$}
This condition will be verified with the help
of the following lemma.

\begin{lemma*} \label{lemma-maximal-closed}
Let $T$ be a scheme and $T'$ a closed subscheme. Let
$\xi':T'\to F$ be a point. Then there is a largest closed
subscheme $\cZ_T\subset T$ such that $\xi'$ extends to $\cZ_T$.
Moreover, its formation is Zariski local: if $U\subset T$
is an open subscheme and $U'=U\cap T'$, we have $\cZ_T\cap U=\cZ_U$.
\end{lemma*}

\begin{proof}
Throughout, for all open subschemes $U\subset T$ we write
$U'=U\cap T'$ and all closed subschemes $Z$ of $U$ such that
$\xi'_{|U'}$ extends to $Z$ are implicitly assumed to contain $U'$.
We proceed by steps.

Let $U=\Spec(A)$ be an affine open subscheme of $T$. Consider the
family of all closed subschemes $Z_\alpha=V(I_\alpha)\subset U$
to which $\xi'_{|U'}$ extends. Consider the ideal $I=\cap I_\alpha$
and define $\cZ_U=V(I)$. Since the map $A/I\to \prod A/I_\alpha$ is
injective, applying Lemma~\ref{lemma:descent-by-sch-dominant}, we
see that $\xi'_{|U'}$ extends to $\cZ_U$. By its very definition
the closed subscheme $\cZ_U$ is largest.

Let $U,V$ be two affine opens of $T$ with $U\subset V$. We claim
that $\cZ_V\cap U=\cZ_U$. Indeed, since $\xi'_{|V'}$ extends to $\cZ_V$
then $\xi'_{|U'}$ extends to $\cZ_V\cap U$, hence $\cZ_V\cap U\subset \cZ_U$.
Conversely, let $Z$ be the schematic image of $\cZ_U\to U\to V$.
The latter map being quasi-compact, the map $\cZ_U\to Z$ is
schematically dominant. By Lemma~\ref{lemma:descent-by-sch-dominant}
it follows that $\xi'_{Z_U}$ extends to $Z$. By maximality this
forces $Z\subset \cZ_V$, hence $\cZ_U\subset \cZ_V\cap U$.

Let $U,V$ be arbitrary affine opens of $T$. We claim that
$\cZ_U\cap V=\cZ_V\cap U$. Indeed, by the previous step, for all affine
opens
$W\subset U\cap V$ we have $\cZ_U\cap V\cap W=\cZ_W=\cZ_V\cap U\cap W$.

Let $\cZ_T$ be the closed subscheme of $T$ obtained by glueing
the $\cZ_U$ when $U$ varies over all affine opens; thus $\cZ_T\cap U=\cZ_U$
by construction. Now $\xi'$ extends to $\cZ_T$, because $\xi'_{|U}$
extends to $\cZ_U$ for each $U$, and we can glue these extensions.
Moreover $\cZ_T$ is maximal with this property, because if $\xi'$
extends to some closed subscheme $Z\subset T$ then for each affine
$U$ the element $\xi'_{|U}$ extends to $Z\cap U$, hence
$Z\cap U\subset \cZ_U=\cZ_T\cap U$, hence $Z\subset \cZ_T$.

The fact that $\cZ_T\cap U=\cZ_U$ for all open subschemes $U$ follows
by restricting to affine opens.
\end{proof}

In order to now verify Condition $(F_8)$, we write $T=\Spec(A)$
and $T'=\Spec(A/I)$. The assumption that
$\xi'_t:\Spec((A/I)_t)\to \Spec(A/I)\to F$ does not lift to any
subscheme of $\Spec(A_t)$ which is strictly larger than
$\Spec((A/I)_t)$ means that the inclusion $T'\subset\cZ_T$
is an equality at the generic point. It follows that
$T'\cap W=\cZ_T\cap W=\cZ_W$ for some open $W$.
Applying Lemma~\ref{lemma-maximal-closed} to variable opens
$W_1\subset W$, we obtain $T'\cap W_1=\cZ_T\cap W_1=\cZ_{W_1}$,
which shows that $W$ fulfills the required condition and we are done.
\end{notitle*}

\begin{notitle*}{Conclusion of the proof}
\label{ssection:conclusion}
Let $G$ be a general finitely presented $S$-group scheme of
multiplicative type~$G$. Let $G_n=\ker(n:G\to G)$ be the
finite flat torsion subschemes; the limit $L=\lim \Hom(G_n,H)$
is an affine $S$-scheme. By reducing to the case $G=\GG_m$
and using the monomorphism $\Hom(G,H)\to L$, we have thus
proven that $\Hom(G,H)$ is representable by a scheme which
is locally of finite presentation over $L$ and over $S$.

To finish the proof of item~(2) of Theorem~\ref{th-2}, let
$F$ be $\Hom(G,H)$ viewed as an $L$-scheme. Let $Y\subset F$
be an $S$-quasi-compact closed subscheme.
Let $Z\subset L$ be the schematic image of $Y\to L$. This
is a closed subscheme and by quasi-compactness, the morphism
$Y\to Z$ is schematically dominant. It follows from
Lemma~\ref{lemma:descent-by-sch-dominant} that $F(Z)\to F(Y)$
is bijective, that is $Y\to F$ factors uniquely through
a map $Z\to F$. Since $Z\into L$ is a closed immersion and
$F\to L$ is separated, then $Z\to F$ is a closed immersion.
By the same argument $Y\to Z$ is a closed immersion;
being also schematically dominant, it is an isomorphism.
Thus~$Y$ is in fact a closed subscheme of $L$. Since $L$
is affine over $S$, it follows that~$Y$ is affine over $S$.
\end{notitle*}

\section{Representability of the functor of subgroups}

In this section we prove item (1) of Theorem~\ref{th-2}.
It would be possible to do this again using Grothendieck's
theorem on unramified functors, with arguments very
similar to those used to prove item (2). However, to
save energy and spare the reader tedious repetitions,
we reduce (1) to (2) by using more advanced technology.
We write $\Sub(H)$ for the functor of multiplicative type
subgroups of $H$.

\begin{notitle*}{Disjoint sum decomposition}
The type $M$ of a group of multiplicative type is
locally constant (\cite{SGA3.2}, Exp.~IX, Rem.~1.4.1).
Although we will note need this, note that by looking at the
inclusion $G_s\subset H_s$ in a fibre, we see that only abelian
groups $M$ of finite type occur. Hence
\[
\Sub(H)=\underset{M}{\textstyle\coprod} \Sub_M(H)
\]
where $\Sub_M(H)$ is the functor of multiplicative type
subgroup schemes of $H$ of type $M$. Thus it is enough to
establish that $\Sub_M(H)$ is representable.
\end{notitle*}

\begin{notitle*}{Representability by an algebraic space}
Let $D(M)=\Spec \cO_S(M)$ be the diagonalizable group
of type $M$ as in \cite{SGA3.2}, Exp.~VIII.
By item (2) of Theorem~\ref{th-2}, the functor
$\Hom(D(M),H)$ is representable by a scheme. It follows from
\cite{SGA3.2}, Exp.~IX, Cor.~6.6 that the subfunctor
\[
\Mono(D(M),H)\subset\Hom(D(M),H)
\]
of monomorphisms of group schemes is an open subscheme.
Moreover, since $D(M)$ is of multiplicative type and $H$
is affine, by \cite{SGA3.2}, Exp.~IX, Cor.~2.5 any
monomorphism $f:D(M)\to H$ is a closed immersion, inducing
an isomorphism between $D(M)$ and a closed subgroup scheme
$K\into H$. By taking a monomorphism $f$ to its image $K$,
we obtain a morphism of functors:
\[
\pi:\Mono(G,H)\to \Sub_M(H).
\]
Let $A=\Aut(D(M))$ be the functor of automorphisms
of $D(M)$; this is isomorphic
to the locally constant group scheme $\Aut(M)_S$.
It acts freely on $\Mono(D(M),H)$ by the rule
$af=f\circ a^{-1}$ for $a\in\Aut(D(M))$ and $f\in\Mono(D(M),H)$.
Let $\Mono(D(M),H)/A$ be the quotient sheaf; this is an
algebraic space by Artin's Theorem, see \cite{SP20},
Tag~\href{https://stacks.math.columbia.edu/tag/04S5}{04S5}.
Since the morphism $\pi$  is $A$-equivariant, it induces
a morphism
\[
i:\Mono(D(M),H)/A\to \Sub_M(H).
\]
We claim that $i$ is an isomorphism. It is enough to
prove that it is an isomorphism of fppf sheaves:
\begin{itemize}
\item surjectivity: this can be checked \'etale-locally
over $S$, so we can assume that the given subgroup scheme
$K\subset H$ is isomorphic to $D(M)$ and then the canonical
inclusion of $D(M)$ into $H$ is a monomorphism that
provides a point of $\Mono(D(M),H)$ lifting $K$.
\item injectivity: if $f_i:D(M)\to H$ are two
monomorphisms with the same image $K$, then
$f_2^{-1}\circ f_1:D(M)\to K\to D(M)$ is an
automorphism of $D(M)$. This settles the claim.
\end{itemize}
\end{notitle*}

\begin{notitle*}{Representability by a scheme}
To prove that $\Sub_M(H)$ is representable by a scheme,
let $L:=\lim \Sub_{M/nM}(H)$ be the limit of the functors of finite
flat multiplicative type subgroups of type $M/nM$. Since
$\Sub_{M/nM}(H)$ is representable and affine (\cite{SGA3.2},
Prop.~3.12.a), the functor $L$ is an affine scheme. By mapping
any subgroup $G\into H$ to the collection of subgroups $G_n\into H$
where $G_n=\ker(n:G\to G)$, we define a morphism of functors
$u:\Sub_M(H)\to L$. By the Density Theorem, this is a
monomorphism. As $\Sub_M(H)\to S$ is locally of finite type,
so is $u$. In particular $u$ is a separated, locally
quasi-finite morphism. By \cite{SP20},
Tag~\href{https://stacks.math.columbia.edu/tag/0418}{0418}
all such morphisms are representable by schemes, hence
$\Sub_M(H)$ is a scheme. Finally, in order to prove that
each $S$-quasi-compact closed subscheme of $\Hom(G,H)$
is affine over $S$, we proceed as in~\ref{ssection:conclusion}.
\end{notitle*}

\bigskip

\noindent
Matthieu ROMAGNY,
{\sc Univ Rennes, CNRS, IRMAR - UMR 6625, F-35000 Rennes, France} \\
Email address: {\tt matthieu.romagny@univ-rennes1.fr}

\end{document}